\documentclass[12pt]{article}
\usepackage[cp1251]{inputenc}
\usepackage[english]{babel}
\usepackage{amsmath,latexsym,amsthm,amsfonts}
\usepackage{amssymb,amscd}
\usepackage{amscd}
\usepackage{caption}
\usepackage{latexsym}

\usepackage[usenames]{color}
\usepackage{colortbl}

\usepackage{graphicx}

\usepackage{longtable}
\usepackage{amssymb}
\usepackage{multicol}
\usepackage{ifthen}
\usepackage{wrapfig}
\usepackage{hyperref}
\usepackage[all]{xy}

\begin{document}

\renewcommand{\figurename}{Fig.}

\begin{center}
{\Large \textbf{Possibility of the use of Cartesian method in the
proofs of fundamental theorems of school planimetry}.}
\end{center}

\begin{center}
Makar Plakhotnyk\\
Post doctoral researcher at S\~ao Paulo University (Brazil)\\
mail:\, Makar.Plakhotnyk@gmail.com
\end{center}

\begin{abstract}
We show how Cartesian method can be used in the proof of
fundamental planimetric topics of the school course, such as
introduction of trigonometric functions, equation of a line and
similarity of triangles.

This work also can be considered as a plan of the school course of
geometry, where the Cartesian method plays the main role.
\end{abstract}

\section*{Introduction}

The work~\cite{Plakh-2016} is published in Ukrainian and concerns
the manner of use of the Thales' theorem in the courses of
geometry in Ukrainian schools.

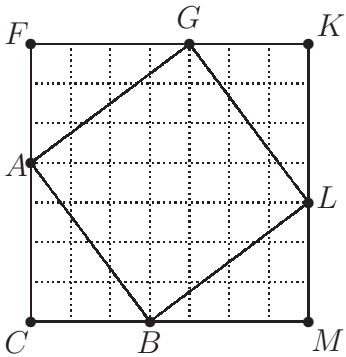
\begin{wrapfigure}{l}{5.5cm}
\setlength{\unitlength}{1pt}
\begin{picture}(125,125)

\put(10,12){\line(1,0){105}} \put(10,12){\circle*{4}}
\put(10,12){\line(0,1){105}} \put(10,117){\line(1,0){105}}
\put(115,117){\line(0,-1){105}}

\put(115,12){\circle*{4}}

\put(55,12){\circle*{4}} \put(10,72){\circle*{4}}
\put(70,117){\circle*{4}} \put(115,57){\circle*{4}}

\put(10,117){\circle*{4}} \put(115,117){\circle*{4}}

\put(0,0){$C$} \put(50,0){$B$} \put(0,67){$A$}

\multiput(10,27)(0,15){6}{\begin{picture}(105,1)
\qbezier[50](0,0)(52.5,0)(105,0)
\end{picture}}
\multiput(25,12)(15,0){6}{\begin{picture}(1,105)
\qbezier[50](0,0)(0,52.5)(0,105)
\end{picture}}

\put(65,123){$G$} \put(0,117){$F$} \put(118,120){$K$}
\put(118,55){$L$} \put(116,1){$M$}

\qbezier(55,12)(55,12)(10,72) \qbezier(10,72)(10,72)(70,117)
\qbezier(70,117)(70,117)(115,57) \qbezier(115,57)(115,57)(55,12)

\end{picture}
\caption{The proof of the Pythagorean theorem ``via
Areas''}\label{fig_2a}
\end{wrapfigure}

Here we present the mathematical ideas of~\cite{Plakh-2016} and
these ideas deal with the possibility and convenience of the use
of Cartesian method in the proof of the Thales' theorem (about the
ratios of various line segments that are created if two
intersecting lines are intercepted by a pair of parallels),
introduction of trigonometric functions and proof the criterions
of similarity of triangles.

It is well known, that Thales' theorem can be proved with the use
of areas. Also is is known thousands proofs of the Pythagorean
theorem. One of these proofs uses the notion of area. Another
classical proof of the Pythagorean theorem uses the notion of the
correctness of definition of $\sin$ and $\cos$ of obtuse angles as
a relation of lengthes of correspond sides of the right triangle.
In the first two chapters of our work we will remind these proofs
for the simplicity of reader.

As our main results we will show the following.

1. It follows from the Pythagorean theorem that the graph of the
equation $y=ax$ for a parameter $a\in\mathbb{R}$ is a line.

2. Suppose that we know that every non-vertical lines on cartesian
plain, passing through origin, is the graph of the equation $y=ax$
for some $a$. This fact directly yields the Thales' theorem.

3. Thales' theorem immediately follows from the Pythagorean
theorem.

This work was motivates by the following facts. In Ukraine the
school course of mathematics is traditionally divided to algebra
and geometry.

In both modern ukrainian books of geometry~\cite{Apostolowa, Bewz,
Burda, Jerszow, Kiseliow, Merzliak-zv, Merzliak}, both in more
ancient Soviet~\cite{Aleksandrow, Kiseliow, Pogorielow}
Pythagorean theorem is proves with the use of similarity of
triangles, which is equivalent to correctness of trigonometric
functions $\sin$ and $\cos$ of obtuse angle and also is equivalent
to the Thales' theorem.

The fact, that the graph of the function $y=ax$ is a line is
studied in ukrainian schools in the course of algebra but this
fact is not properly proved in the books~\cite{Bewz-7, Ister-7,
Krawczuk-7, Merzliak-7}.

\section{Proof of the Pythagorean theorem with the use of areas}

Consider a $ABC$ triangle with a right angle $C$ (see
pict.~\ref{fig_2a}). Take points $F$ and $M$ on the lines $CA$ and
$CB$ such that segments $CF$ and $CM$ be equal to the sum of legs
of the triangle $ABC$.

Take points $G$ and $L$ on the sides on the square $CFKM$ such
that $FG = KL =AB$. Then the four new triangles will be equal, and
the quadrangle $AGLB$ will be a square.

The proof of the Pythagorean theorem can be obtained, if we
subtract the area of the four right triangles from the area of the
``big'' square.

\section{Proof of the Thales' theorem using areas}

The proof of the correctness of trigonometric functions in the
most widespread school books in Ukraine is realized by via the
Generalized Thales' theorem. Thus, we will give below the proof of
the correctness of trigonometric functions, which is given in the
by Aleksandroff in~\cite[p. 95]{Aleksandrow}.

Take any right triangle $ABC$ (see pict~\ref{fig_2c}) and a point
$M$ on the leg $AC$. Express the area of the triangle $ABM$ by
manners. From one hand, $S = \frac{1}{2}ma$, where $a=BC$ and $m =
AM$. From another hand, $S = \frac{1}{2}ch$, where $h = MD$ is the
height of the triangle $ABM$ and $c$ is the hypotenuse of the
former right triangle. Sine the obtain formulas represent the same
area, then $$ \frac{a}{c} = \frac{h}{m}.
$$

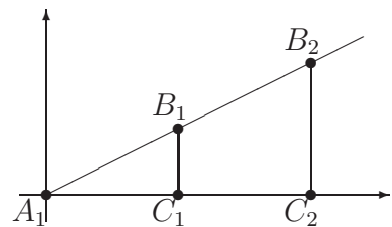
\section{Introduction of trigonometric functions using
cartesian method}

Let triangles $A_1B_1C_1$ and $A_2B_2C_2$ have right angles
$\angle C_1 = \angle C_2 = 90^0$ and let $\angle A_1 = \angle
A_2$.

Superpose points $A_1$ and $A_2$ and assume that $C_1$ and $C_2$
are on the same side from $A_1  = A_2$ on the line $A_1C_1$. Take
the cartesian plain with the origin in $A = A_1$ and the x-axis
$A_1C_1$ such that $B_1$ and $B_2$ belong to the first quadrant
(see fig.~\ref{fig_2b}).

In this case semi-straight lines $AB$ and $A_1B_1$ will coincide.
The equation of their line is $y = ax$, where real $a$ depends on
the value of angle $\angle A = \angle A_1$.

Figure~\ref{fig_2b} is prepared for the case when $A_1C_1<A_2C_2$.
Denote $A_1C_1 = x_1$ and $A_2C_2 = x_2$. Then we can express the
coordinates of the vertices of our triangles as follows: $A_1(0,\,
0)$, $A_2(0,\, 0)$, $B_1(x_1,\, ax_1)$, $B_2(x_2,\, ax_2)$,
$C_1(x_1,\, 0)$ and $C_2(x_2,\, 0)$.

Obtain the lengthes of the sides of triangles from the distance
formula, precisely $A_1C_1 = x_1$, $A_2C_2 = x_2$, $B_1C_1 =
ax_1$, $B_2C_2 = ax_2$, $A_1B_1 = x_1\sqrt{1+a^2}$ and $A_2B_2 =
x_2\sqrt{1+a^2}$.

\begin{figure}[htbp]
\begin{minipage}[h]{0.45\linewidth}
\begin{center}
\begin{picture}(200,140)

\put(10,15){\line(1,0){200}} \put(10,15){\line(3,2){200}}
\put(190,15){\line(0,1){120}} \put(110,15){\line(-2,3){30.5}}
\qbezier(110,15)(110,15)(190,135)
\qbezier(150,108)(150,108)(150,15) \put(150,108){\circle*{4}}
\put(150,15){\circle*{4}}

\put(110,15){\circle*{4}} \put(10,15){\circle*{4}}
\put(79.5,61){\circle*{4}} \put(190,135){\circle*{4}}
\put(190,15){\circle*{4}}

\put(3,2){$A$} \put(100,2){$M$} \put(180,2){$C$} \put(68,63){$D$}
\put(180,138){$B$}

\put(193,75){$a$} \put(60,18){$m$} \put(95,40){$h$}
\put(153,55){$a'$} \put(140,2){$C'$} \put(140,110){$B'$}
\end{picture}\end{center}
\caption{The correctness of trigonometric functions by
Alexandroff}\label{fig_2c}
\end{minipage}
\hfill
\begin{minipage}[h]{0.45\linewidth}
\begin{center}
\begin{picture}(200,140)

\put(0,15){\vector(1,0){140}}

\put(10,15){\circle*{4}} \put(-3,5){$A_1$}

\put(10,5){\vector(0,1){80}}

\put(10,15){\line(2,1){120}} \put(110,65){\line(0,-1){50}}

\put(60,40){\line(0,-1){25}}

\put(60,40){\circle*{4}} \put(50,45){$B_1$}
\put(110,65){\circle*{4}} \put(100,70){$B_2$}

\put(60,15){\circle*{4}} \put(50,5){$C_1$}

\put(110,15){\circle*{4}} \put(100,5){$C_2$}

\end{picture}
\end{center} \caption{Introduction of trigonometric functions}\label{fig_2b}

\end{minipage}
\hfill
\end{figure}

Due to these formulas, the equalities of ratios of the lengthes of
sides of the triangles, which are necessary for the introduction
of trigonometric functions, are evident.

\section{Similarity of right triangles}

When the cartesian plain and equation of a line are introduced,
the criterions of similarity of triangles become obvious facts,
which follow from Cartesian method.

Remind that triangles are called similar, if their correspond
sides are proportional and correspond angles are equal. Remind the
criterions of similarity of right triangles.

\begin{enumerate}

\item By two angles (two right triangles are similar if and only
if a pair of their correspond angles are equal).

\item By two pairs of proportional sides (two right triangles are
similar if and only if a pair of their correspond sides have the
proportional lengthes).
\end{enumerate}

The first of these criterions follows from reasonings, which are
analogous the one from the introducing of trigonometric functions.
It is necessary to consider the Cartesian plane, superpose
vertices of two equal acute angles at origin, plug right angles to
the x-axis and plug hypotenuses at the same line, passing through
origin. Not the fact that corresponds sides of triangles are
proportional follows from the distance formula and the formula of
equation of the line, passing through origin.

The second criterion of the similarity of right triangles can be
proved in the same manner. Suppose that sides of triangles
$A_1B_1C_1$ are $A_2B_2C_2$ with right angles $\angle C_1$ and
$\angle C_2$ are proportional. Like we have done above, superpose
vertices $A_1$ and $A_2$ at origin, plug $C_1$ and $C_2$ at x-axis
and plug $B_1$ and $B_2$ in the first quadrant. Suppose the
equation of lines $A_1B_1$ and $A_2B_2$ be $y=ax$ and $y=bx$
respectively. Nevertheless, $a = \frac{B_1C_1}{A_1C_1} =
\frac{B_2C_2}{A_2C_2} = b,$ which implies the equality of all
correspond angles of our triangles.

\section{Similarity of triangles in general case}

We suggest to reduce the ``general case'' of the question on
similarity of triangles to the similarity of right triangles.

We will introduce the idea of the proof of the similarity of
triangles in ``general case''.

\textbf{Similarity by three angles}. Suppose that correspond
angles of triangles $A_1B_1C_1$ and $A_2B_2C_2$ are equal. Then
sides of triangles are proportional.

Consider the height of the biggest angle of these triangle, let it
be $\angle C_1 = \angle C_2$. Since the angle if the biggest, then
the base of the height will belong to the side of triangle, but
not to its continuation.

It follows from the similarity of right triangles by three angles,
that $\triangle B_1C_1H_1 \sim \triangle B_2C_2H_2$ and $\triangle
A_1C_1H_1 \sim \triangle A_2C_2H_2$, whence $$
\underbrace{\frac{a_1}{a_2} = \frac{x_1}{x_2} =
\frac{h_1}{h_2}}_{\triangle B_1C_1H_1 \sim \triangle B_2C_2H_2} =
\frac{y_1}{y_2} = \frac{b_1}{b_2}.
$$ Thus, $x_1 = \frac{h_1 x_2}{h_2}$ and $y_1 =
\frac{h_1y_2}{h_2}$, which implies that $\frac{x_1+y_1}{x_2+y_2}
=\frac{h_1}{h_2}$ and correspond sides of $A_1B_1C_1$ are
$A_2B_2C_2$ proportional.

\textbf{Similarity by an angle and two proportional sides}.
Suppose that triangles $A_1B_1C_1$ and $A_2B_2C_2$ are such that
$\angle B_1 = \angle B_2$ and sides of the angle $B_1$ are
proportional to side of~$B_2$. Consider the heights from vertices
$C_1$ and $C_2$ to the correspond side, or its continuation.

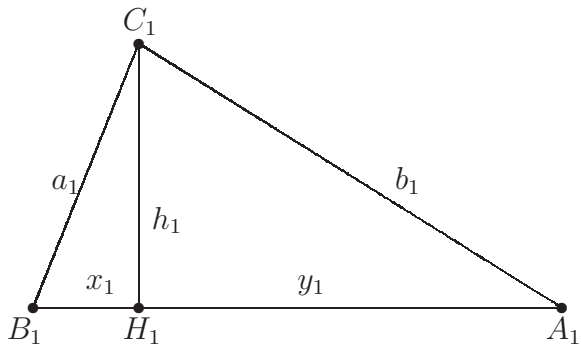
\begin{wrapfigure}{l}{8cm}
\begin{center}
\setlength{\unitlength}{1pt}
\begin{picture}(200,125)

\put(0,0){\line(1,0){200}} \qbezier(0,0)(0,0)(40,100)
\qbezier(40,100)(40,100)(200,0) \put(40,100){\line(0,-1){100}}

\put(-10,-12){$B_1$} \put(194,-12){$A_1$} \put(34,-12){$H_1$}
\put(34,105){$C_1$}

\put(0,0){\circle*{4}} \put(40,100){\circle*{4}}
\put(40,0){\circle*{4}} \put(200,0){\circle*{4}} \put(7,45){$a_1$}
\put(137,45){$b_1$} \put(20,7){$x_1$} \put(100,7){$y_1$}
\put(45,30){$h_1$}

\end{picture}
\end{center} \caption{On the
similarity of triangles}\label{fig_1-podibn-kuty}
\end{wrapfigure}

Similarity $\triangle B_1C_1H_1 \sim \triangle B_2C_2H_2$ follows
from the criterion ``by three angles'. This, together with
proportionality of sides of angles $B_1$ and $B_2$ imply that
$$ \left\{
\begin{array}{l}
\frac{h_1}{h_2} = \frac{x_1}{x_2} = \frac{a_1}{a_2},\\
\frac{a_1}{a_1} = \frac{x_1+y_1}{x_2+y_2}.
\end{array}\right.
$$

It follows from these equalities that legs of right triangles
$C_1H_1A_1$ and $C_2H_2A_2$ are proportional, whence triangles are
similar.

Similarity of right triangles imply that correspond angles of
$A_1B_1C_1$ and $A_2B_2C_2$ are equal, whence triangles are
similar.

\textbf{Similarity by proportional sides}. Suppose that sides of
triangles $A_1B_1C_1$ and $A_2B_2C_2$ are proportional and prove
the equality of angles.

Consider the triangle
$\widetilde{A}_2\widetilde{B}_2\widetilde{C}_2$, which is similar
to $A_1B_1C_1$ and equal to $A_2B_2C_2$. For this deal take an
arbitrary points  $\widetilde{A}_2$ and $\widetilde{B}_2$ such
that $\widetilde{A}_2\widetilde{B}_2 = A_2B_2$ and take rays to
the same semi plain of $\widetilde{A}_2\widetilde{B}_2$, with
angles $\angle A_1$ and $\angle B_1$ between rays and the line.
Denote by~$\widetilde{C}_2$ the point of intersection of these
rays.

Triangle $\widetilde{A}_2\widetilde{B}_2\widetilde{C}_2$ is
similar to $A_1B_1C_1$ by three angles, whence sides of triangles
are proportional. From another hand, triangles
$\widetilde{A}_2\widetilde{B}_2\widetilde{C}_2$ and $A_2B_2C_2$
are equal by three sides, since $\widetilde{A}_2\widetilde{B}_2 =
A_2B_2$ and correspond sides of these triangles are proportional
to sides of $A_1B_1C_1$.

\section{Proof of the Thales' theorem
as a corollary of Pythagorean theorem}

We will show in this section how Thales' theorem implies from the
Pythagorean theorem.

\begin{figure}[htbp]
\begin{minipage}[h]{0.45\linewidth}
\begin{center}
\begin{picture}(200,230)

\put(10,0){\line(1,0){100}} \put(10,0){\line(0,1){200}}
\put(10,0){\circle*{4}} \put(10,200){\circle*{4}}
\put(110,0){\circle*{4}}

\put(50,0){\line(0,1){120}} \put(50,0){\circle*{4}}
\put(50,120){\circle*{4}} \put(10,200){\line(1,-2){100}}

\put(10,120){\line(1,0){40}} \put(10,120){\circle*{4}}

\put(28,3){$a$} \put(75,3){$b$} \put(15,153){$c$} \put(15,50){$d$}
\put(55,50){$d$} \put(28,123){$a$}

\put(0,-13){$A$}  \put(45,-13){$D$} \put(0,202){$B$}
\put(105,-13){$C$} \put(55,122){$F$} \put(0,122){$H$}

\end{picture}
\end{center} \caption{The special case of the Thales' theorem}
\label{fig_3a}
\end{minipage}
\hfill
\begin{minipage}[h]{0.45\linewidth}
\begin{center}
\begin{picture}(200,230)

\put(0,0){\line(1,0){200}} \put(0,0){\line(2,5){40}}

\qbezier(40,100)(40,100)(200,0)

\put(40,100){\line(0,-1){100}}

\put(20,50){\line(1,0){100}}

\put(-10,-13){$B$} \put(194,-13){$A$} \put(34,-13){$H$}
\put(34,105){$C$}

\put(0,0){\circle*{4}} \put(40,100){\circle*{4}}
\put(20,50){\circle*{4}} \put(40,50){\circle*{4}}
\put(120,50){\circle*{4}} \put(40,0){\circle*{4}}
\put(200,0){\circle*{4}} \put(0,30){$x_1$} \put(160,30){$y_1$}
\put(80,80){$y_2$} \put(20,80){$x_2$} \put(43,20){$h_1$}
\put(43,75){$h_2$} \put(13,55){$F$}
\put(43,38){$G$}\put(125,55){$K$}
\end{picture}
\end{center} \caption{The proof of Thales'
theorem in the general case}\label{fig_3b}
\end{minipage}
\hfill
\end{figure}

Use the figure~\ref{fig_3a}, where small letters $a,\, b,\, c,\,
d$ denote segments $AD=HF$, $DC$, $HB$ and $AH = DF$. Start from
the case when lines  $BA$ and $FD$, which intersect sides of the
angle $\angle BCA$, are perpendicular to $AC$.

Use Pythagorean theorem to express $FC$ from the right triangle
$DFC$, $FB$ from $HFB$ and $BC$ from $ABC$. Then equality
$$ FC + FB = BC$$ can be rewritten as $$ \sqrt{d^2 + b^2} +
\sqrt{a^2 + c^2} = \sqrt{(a+b)^2 + (d+c)^2}.
$$

After evident transformations obtain $(ad - bc)^2 =0. $ Use
Pythagorean theorem ones more, rewrite this equality as $$
\frac{\sqrt{a^2+c^2}}{a} = \frac{\sqrt{b^2+d^2}}{b},
$$ which is one of the forms of the Thales' theorem.

The proof of Thales' theorem without the assumption that lines are
perpendicular to one of the side of an angle is the following. Use
the notation of the figure~\ref{fig_3b}, which is prepared for the
case when perpendicular from the vertex of the angle to parallel
lines is inside, but not outside of the angle. For another case
the proof is almost the same.

Apply Thales' theorem for the angles $BCH$ and $HCA$ obtain
equalities $$ \frac{x_2}{x_1} = \frac{h_2}{h_1} = \frac{y_2}{y_1},
$$ which imply Tales theorem for general case.

\end{document}